\documentclass[a4paper, 12pt]{article}
\usepackage{cmap}					
\usepackage[T2A]{fontenc}			
\usepackage[utf8]{inputenc}			
\usepackage[english]{babel}
\usepackage{amsthm,amsmath,amsfonts,amssymb}
\usepackage{color}
\usepackage{fontenc}

\usepackage[left=2.25cm,right=2.25cm,top=2.5cm,bottom=2.5cm,bindingoffset=0cm,nohead,nofoot,footskip=7mm]{geometry}

\newtheorem{lemma}{{\bf Lemma}}[section]
\newtheorem{theorem}{{\bf Theorem}}[section]
\newtheorem{rem}{{\bf Remark}}

\author{
A.I. Nazarov\footnote{St. Petersburg Dept of Steklov Mathematical Institute; St. Petersburg State University. E-mail: al.il.nazarov@gmail.com }\,,
A.P. Scheglova\footnote{St. Petersburg Electrotechnical University. E-mail: alexandra.scheglova@gmail.com}}

\title{On the sharp constant\\ in ``magnetic'' 1D embedding theorem }

\date{}

\begin{document}

\maketitle

\section{Introduction}

We consider the problem of finding the sharp (exact) constant in the ``magnetic'' embedding theorem 
\begin{equation}\label{f_001}
\min\limits_u\frac {\Vert u'+iA u\Vert_{L_2}}{\Vert u\Vert_{L_q}}=: \mu_q (A),
\end{equation}
where $A\in L_1(0,2\pi)$, and minimum is taken over all $2\pi$-periodic absolutely continuous functions. 

It is easy to see that $\mu_q(A)$ is attained and does not change if we change $A\mapsto A+k$, $k\in\mathbb Z$. Moreover, the substitution
$$
u(x)\mapsto u(x)\exp\Big(i\int\limits_0^x (A(t)-\alpha)\,dt\Big), \qquad \alpha=\frac 1{2\pi}\int\limits_0^{2\pi} A(t)\,dt,
$$
shows that we can assume without loss of generality $A\equiv\alpha$ and $|\alpha|\le\frac 12$.\medskip

Trivially the value $\mu_q(0)\equiv0$ is attained by any constant function. Further, if $q\le 2$ then due to the evident estimate 
$\Vert u\Vert_{L_q}\le(2\pi)^{\frac 1q-\frac 12}\cdot\Vert u\Vert_{L_2}$ the constant function also is a minimizer of $\mu_q(\alpha)$, and 
$\mu_q(\alpha)= (2\pi)^{\frac 12-\frac 1q}\cdot|\alpha|$. Thus, the constant function is a natural
candidate to the minimizers of $\mu_q(\alpha)$. In this paper we show that in fact for $\alpha\ne0$ it is minimizer only for sufficiently small $q>2$, namely, for 
$(q+2)\alpha^2\le1$. In particular, for $\alpha=\pm \frac 12$ and $q>2$ the minimizer is always non-constant.

\begin{rem}
For $q=\infty$ the sharp constant in (\ref{f_001}) was found in \cite{GO}, see also \cite{ILLZ}.
\end{rem}

In what follows we assume $2<q<\infty$.
It is convenient to normalize $u$ by 
$\Vert u\Vert_{L_q}^q=2\pi$, and we arrive at the problem
\begin{equation}\label{isoper}
\mu_q(\alpha)^2=(2\pi)^{-\frac 2q}\cdot\min\limits_u 
\int\limits_0^{2\pi}|u'+i\alpha u|^2\,dx,
\qquad 
\int\limits_0^{2\pi}|u|^q\, dx=2\pi.
\end{equation}

To study the problem (\ref{isoper}) we use the phase plane method. In a similar way in \cite{N99}, \cite{N04} the problem
$$
\min\limits_u \int\limits_{-T}^T(u'{}^2+u^2)\,dx, \qquad \int\limits_{-T}^T|u|^q\, dx=1.
$$
was studied, and the sharp condition of symmetry breaking in this problem was found. See also \cite[Lemma 5]{DEL}.\medskip

\section{The constant and non-constant minimizers of (\ref{isoper})}

Denote $r=|u|$ and $\varphi=\arg(u)+\alpha x$. Then (\ref{isoper}) can be rewritten as follows:
\begin{equation}\label{f_003}
J(r,\varphi)=\int\limits_0^{2\pi}|r'+ir\varphi'|^2\,dx=\int\limits_0^{2\pi}(r'{}^2+r^2\varphi'{}^2)\,dx\,\to\,\min,\qquad
\int\limits_0^{2\pi}r^q \,dx=2\pi.
\end{equation}
Here $r$ and $\varphi'$ are $2\pi$-periodic functions, and
\begin{equation}\label{fi}
\int\limits_0^{2\pi}\varphi'\,dx=2\pi\alpha.
\end{equation}

The Euler equation with respect to $\varphi$ reads:
$$
0\equiv\dfrac{1}{2}D_{\varphi}J(r,\varphi)(\psi)=r^2\varphi'\psi\Bigr|_0^{2\pi}-\int\limits_0^{2\pi}\bigl(r^2\varphi'\bigr)'\psi\,dx.
$$
The first term vanishes due to 2$\pi$-periodicity, and we obtain
\begin{equation}\label{f_004}
r^2\varphi'=a=const.
\end{equation}

The Euler--Lagrange equation with respect to $r$ reads:
$$
r'h\Bigr|_0^{2\pi}+\int\limits_0^{2\pi}\bigl(r\varphi'{}^2-\lambda r^{q-1}-r''\bigr)h\,dx\equiv 0,
$$
and we obtain
$$
-r''+r\varphi'{}^2=\lambda r^{q-1}.
$$
Taking into account~(\ref{f_004}) we arrive at
\begin{equation}\label{f_005}
-r''+\frac{a^2}{r^3}=\lambda r^{q-1}.
\end{equation}
It is easy to see that the function $r\equiv1$ is a solution of (\ref{f_005}). Moreover, in this case relations (\ref{f_004}) and (\ref{fi}) give 
$a=\alpha$, and thus $\lambda=\alpha^2$.

\begin{theorem}
Let $(q+2)\alpha^2>1$. Then the function $r\equiv1$ cannot provide minimal value in the problem (\ref{f_003}), and thus we have 
$\mu_q(\alpha)< (2\pi)^{\frac 12-\frac 1q}\cdot|\alpha|$.
\end{theorem}
\noindent{\bf Proof.} Taking into account (\ref{f_004}) we conclude that the second order necessary condition of minimum is positivity of the quadratic form
$$
\int\limits_0^{2\pi}\Bigl(h'{}^2-\dfrac{3 a^2 h^2}{r^4}-\lambda (q-1) r^{q-2}h^2\Bigr)\,dx
$$
on the space of $2\pi$-periodic function with zero mean value. Substituting $r\equiv1$, $a=\alpha$, and $\lambda=\alpha^2$ we obtain
$$
\int\limits_0^{2\pi}\Bigl(h'{}^2-\alpha^2(q+2)h^2\Bigr)\,dx\ge0.
$$
For $(q+2)\alpha^2>1$ this inequality fails for $h=\sin(x)$. \hfill$\blacksquare$\medskip

\begin{theorem}
Let $(q+2)\alpha^2\le 1$. Then the function $r\equiv1$ provides minimal value in the problem (\ref{f_003}), and thus we have 
$\mu_q(\alpha)= (2\pi)^{\frac 12-\frac 1q}\cdot|\alpha|$.
\end{theorem}
\noindent{\bf Proof.} Integrating ODE (\ref{f_005}) we obtain
\begin{equation}\label{f_007}
\dfrac{{r'}^2}{2}=-\dfrac{a^2}{2r^2}-\dfrac{\lambda}{q}r^q+c.
\end{equation}
On the another hand, we can multiply (\ref{f_005}) by $r$ and integrate over $[0,2\pi]$. This gives in view of the normalization condition
$$
\int\limits_0^{2\pi}\left({r'}^2+\dfrac{a^2}{r^2}\right)\,dx=\lambda\int\limits_0^{2\pi}r^q \,dx=2\pi\lambda,
$$
and (\ref{f_007}) implies $c=\frac{\lambda}{2}+\frac{\lambda}{q}$.

If $r$ is not a constant then the right-hand side of (\ref{f_007}) has two zeros corresponding to minimal and maximal values of $r$ at the period. 
Denote these values by $r_1$ and $r_2$ respectively. By the normalization condition we have
\begin{equation}\label{r12}
r_1<1<r_2.
\end{equation}

Thus, any non-constant periodic positive solution of ODE (\ref{f_005}) corresponds to the motion along an oval given by equation (\ref{f_007}) 
in the phase plane $(r,r')$. Since this oval is symmetric w.r.t. $r'$ axis, without loss of generality we can assume that $r(0)=r(2\pi)=r_1$ and $r(\pi)=r_2$.

Consider a half of the oval corresponding to $r'>0$. Then we have 
from (\ref{f_007})
$$
r'=\sqrt{2c-\dfrac{a^2}{r^2}-\dfrac{2\lambda}{q}r^q}=\dfrac{\sqrt{\lambda(1+\frac{2}{q})r^2-a^2-\frac{2\lambda}{q}r^{q+2}}}{r}.
$$
By (\ref{fi}) and (\ref{f_004}) we obtain
\begin{equation}\label{f_008}
2\pi\alpha=a\int\limits_0^{2\pi}\dfrac{dx}{r^2}=2a\int\limits_0^{\pi}\dfrac{dx}{r^2}
=\int\limits_{r_1}^{r_2}\dfrac{2dr}{r\sqrt{\frac{\lambda}{a^2}\left[\left(1+\frac{2}{q}\right)r^2-\frac{2}{q}r^{q+2}\right]-1}}.
\end{equation}

By (\ref{r12}) we have $\frac{\lambda}{a^2}>1$. Changing the variable $t=\frac{\lambda}{a^2}\left(1+\frac{2}{q}\right) r^2$ we rewrite (\ref{f_008}) as follows:
\begin{equation}\label{f_009}
M_q(\gamma):=\int\limits_{t_1}^{t_2}\dfrac{dt}{t\sqrt{t-\gamma t^{\frac q2 +1}-1}}=2\pi\alpha.
\end{equation}
Here $t_1$, $t_2$ are the roots of the equation $t-\gamma t^{\frac q2 +1}-1=0$, 
and 
$$
0<\gamma
<\gamma_{max}=\dfrac{2}{q+2}\left(1+\dfrac{2}{q}\right)^{-\frac{q}{2}}.
$$

The statement of Theorem follows from Lemma which will be proved in Section 3.\medskip

\begin{lemma}\label{lem2}
For all $\gamma\in (0;\gamma_{max})$
we have
\begin{equation}\label{monot}
M_q'(\gamma)<0.
\end{equation}
Moreover,
\begin{equation}\label{ner}
\lim\limits_{\gamma\uparrow \gamma_{max}}M_q(\gamma)=\dfrac{2\pi}{\sqrt{q+2}}.
\end{equation}
\end{lemma}
\medskip

Namely, it follows from (\ref{monot}) and (\ref{ner}) that if $(q+2)\alpha^2\le 1$ then $M_q(\gamma)>2\pi\alpha$ for all
$\gamma\in(0;\gamma_{max})$. Therefore, the equation (\ref{f_009}) has no solutions, and the constant function is a unique stationary point of the problem (\ref{f_003}). 
This completes the proof. \hfill$\blacksquare$

\begin{rem}
 If $(q+2)\alpha^2> 1$ then the equation (\ref{f_009}) has a unique solution. Evidently, the motion along corresponding oval in the phase plane just provides the minimum in 
 (\ref{f_003}).
\end{rem}

\section{Proof of Lemma \ref{lem2}}

We introduce the notation
\begin{equation}\label{f_011}
f(t)=t-\gamma t^{\frac q2 +1}-1.
\end{equation}
Then
\begin{equation}\label{f_014}
f'(t)=1-\dfrac{\gamma(q+2)}{2}\,t^{\frac{q}{2}};\quad
f''(t)=-\dfrac{\gamma q(q+2)}{4}\,t^{\frac{q}{2}-1};\quad
f'''(t)=-\dfrac{\gamma q(q+2)(q-2)}{8}\,t^{\frac{q}{2}-2}.
\end{equation}
It is easy to see that $f'(t_1)>0$ and $f'(t_2)<0$. Denote by $t_0$ a unique root of $f'$. \medskip 

To prove (\ref{ner}) we observe that by the Rolle Theorem for any $t \in\ (t_1,t_2)$ there exists $\overline t(t) \in (t_1,t_2)$ such that
$$
f(t) = - \frac {f''(\overline t)}2\, (t_2-t)(t-t_1).
$$
Hence
$$M_q(\gamma) = \sqrt{\frac {-2}{f''(\widehat t\,)}}\,\int\limits_{t_1}^{t_2}\frac
{dt}{t\sqrt{(t_2-t)(t-t_1)}} =\frac {\pi}{\widetilde t} \sqrt{\frac {-2}{f''(\widehat t\,)}},
$$
where $\widehat t$ and $\widetilde t$ are some points in $(t_1,t_2)$.

Notice that $t_1\uparrow \frac{q+2}{q}$ and $t_2\downarrow \frac{q+2}{q}$ as $\gamma\uparrow\gamma_{max}$. Therefore, 
$\widehat t$ and $\widetilde t$ also tend to $\frac{q+2}{q}$, and
$$
\lim\limits_{\gamma\uparrow \gamma_{max}}M_q(\gamma)=
\frac {\pi q}{q+2} \sqrt{\frac {-2}{f''(\frac{q+2}{q})\big|_{\gamma=\gamma_{max}}}}=
\frac{2\pi}{\sqrt{q+2}},
$$ 
and (\ref{ner}) follows.\medskip

To prove (\ref{monot}) we proceed similarly to \cite[Sec. 2]{N99} and \cite{N04}.

\begin{lemma} 
For $ \gamma \in (0, \gamma_{max})$ the following identity holds:
\begin{equation}\label{M'}
M'_q(\gamma) = \int\limits_{t_1}^{t_2}\frac{\sqrt{f} f'}{\Psi^2}\cdot H_{\beta}\, dt, 
\end{equation}
where 
\begin{eqnarray*}
\Psi&=&f'^2-2ff'',\\
H_{\beta}&=&
\beta(3{f'}^2f''+2ff'f'''-6f{f''}^2)-\frac{q(q-2)t^{\frac{q}{2}-2}}{2},
\end{eqnarray*}
and $\beta$ is an arbitrary number.
\end{lemma}

\noindent{\bf Proof.} 
We have
$$ M_q^{(\epsilon)} (\gamma) :=\int\limits_{t_1+\epsilon}^{t_2-\epsilon} \frac{dt}{t\sqrt{f}}\ \stackrel
{\epsilon\to0} \longrightarrow \ M_q(\gamma),
$$
and convergence is uniform in any compact subset of the interval $(0,\gamma_{max})$.

Furthermore,
$$\frac{dM_q^{(\epsilon)}}{d\gamma} = \frac{dt_2}{d\gamma}\cdot \left.\frac 1{t\sqrt{f}}\,
\right\vert^{t_2-\epsilon} -\ \frac{dt_1}{d\gamma}\cdot\left.\frac 1{t\sqrt{f}}\, \right\vert^{t_1+\epsilon} +\ 
\frac 12\int\limits_{t_1+\epsilon} ^{t_2-\epsilon} \frac {t^{\frac q2} \, dt}{f^{\frac 32}}.
$$

However, $f(t_1)=f(t_2)=0$ implies
$$
 \partial_{\gamma}f\big|^{t_k}+f'\big|^{t_k}\cdot\frac{dt_k}{d\gamma}=0,\quad k=1,2.
$$
Therefore,
$$\frac{dt_k}{d\gamma}=\left.\frac{t^{\frac q2 +1}}{f'}\,\right\vert^{t_k}=
\left.\frac{t^{\frac q2 +1} f' - qt^{\frac q2}f+2\beta f^2f''}{f'^2-2ff''}\,\right\vert^{t_k},\quad k=1,2,
$$
and thus
$$\frac{dM_q^{(\epsilon)}}{d\gamma} = \left.\frac 1{t\sqrt{f}}\ \frac{t^{\frac q2}(tf'-qf)+2\beta f^2f''}{\Psi}\,
\right\vert_{t_1+\epsilon}^{t_2-\epsilon}\ +\ O(\epsilon^{\frac 12})\ +\ 
\frac 12\int\limits_{t_1+\epsilon}^{t_2-\epsilon} \frac {t^{\frac q2}\, dt}{f^{\frac 32}}.
$$

Note that $\Psi(t_1)=f'^2(t_1)>0$, and
$$\Psi'= -2f\cdot f'''>0\qquad \mbox{in}\quad  (t_1,t_2).
$$
Hence $\Psi >0$ in $ [t_1,t_2]$, and we can write
$$
\frac{dM_q^{(\epsilon)}}{d\gamma}=\int\limits_{t_1+\epsilon}^{t_2-\epsilon}\,\left[\frac {d}{dt} 
\left(\dfrac{1}{t\sqrt{f}}\cdot \dfrac{t^{\frac q2}(tf'-qf)+2\beta f^2f''}{\Psi}\right)+\frac{t^{\frac q2}}{2f^{\frac 32}}\,\right]dt\ +\ O(\epsilon^{\frac 12}).
$$
The expression in square brackets is equal to $\dfrac{\sqrt{f} f'}{\Psi^2}\cdot H_{\beta}$.
Therefore $dM_q^{(\epsilon)}/d\gamma$ converges to the right-hand side of (\ref{M'}) as $\epsilon \to 0$. 
Moreover, convergence is uniform  in any compact subset of the interval $(0,\gamma_{max})\,$. 
This completes the proof.\hfill$\blacksquare$\medskip

Using the relations (\ref{f_011})--(\ref{f_014}) we calculate
$$
H_{\beta}(t)= q t^{\frac{q}{2}-2}
\Bigl(\dfrac{\beta \gamma (q+2)}{16} h(t)-\dfrac{q-2}{2}\Bigr),
$$
where
$$
h(t)=4(q-2)-4(q+1)t-4\gamma(q+1)(q-2)t^{\frac{q}{2}+1}+4\gamma(q+2)(q+1)t^{\frac{q}{2}}+\gamma^2(q+2)(q-2)t^{q+1}.
$$

Direct calculation shows that
$$
h''(t)= 
-\gamma q(q+1)(q-2)(q+2)t^{\frac{q}{2}-2}f(t). 
$$
Thus, $h''(t)<0$ for $t\in(t_1,t_2)$. Therefore, $h'$ decreases on $[t_1,t_2]$.

Next, the relation (\ref{f_011}) implies $\gamma t_1^{\frac{q}{2}}=\frac{t_1-1}{t_1}$. Therefore,
\begin{eqnarray*}
h'(t_1)&=&-(q+1)\Bigl(4+2\gamma(q-2)(q+2)t_1^{\frac{q}{2}}-2\gamma q(q+2)t_1^{\frac{q}{2}-1}-\gamma^2(q+2)(q-2)t_1^{q}\Bigr)\\
&=&-(q+1)\biggl[4+(q+2)\biggl(2(q-2)\frac{t_1-1}{t_1}-2q\frac{t_1-1}{t_1^2}-(q-2)\Bigl(\frac{t_1-1}{t_1}\Bigr)^2\biggr)\biggr]\\
&=&
-\frac{q+1}{t_1^2}\bigl(qt_1-(q+2)\bigr)^2<0.
\end{eqnarray*}
Thus, $h'(t)<0$ on $[t_1,t_2]$. It follows that for $\beta<0$ the function $H_{\beta}$ increases on $[t_1,t_2]$.

Now we choose
$$
\beta=-\frac{q(q-2)t_0^{\frac{q}{2}-2}}{12f(t_0){f''}^2(t_0)}<0
$$
(we recall that $f'(t_0)=0$). Then  $H_{\beta}(t_0)=0$. By monotonicity we have $H_{\beta}<0$ on $[t_1;t_0)$ and $H_{\beta}>0$ on $(t_0;t_2]$. Therefore,
$\dfrac{\sqrt{f}f'}{\Psi^2}\cdot H_{\beta}\le0$ on $[t_1,t_2]$, and (\ref{M'}) implies (\ref{monot}). \hfill$\blacksquare$

\section*{Acknowledgements}

We are grateful to Prof. Ari Laptev for the statement of the problem and useful discussion. A.N. thanks the Mittag-Leffler Institute for the hospitality during
the visit in January 2017.

Authors' work was supported by RFBR grant 17-01-00678a.

\end{document}